\documentclass{birkjour}
\usepackage{amsmath}
\usepackage{amssymb}
\usepackage{amsthm}
\usepackage{color}
\usepackage{eqlist}
\usepackage{enumerate}
\usepackage{txfonts}
\usepackage{geometry}
\usepackage[mathscr]{eucal}
\usepackage{hyperref}
\hypersetup{
     colorlinks   = true,
     citecolor    = blue
}
\hypersetup{linkcolor=blue}

 \newtheorem{thm}{Theorem}[section]
 
 \newtheorem{lem}[thm]{Lemma}
 \newtheorem{prop}[thm]{Proposition}
 \theoremstyle{definition}
 
  \newtheorem{ex}[thm]{Example}
 \newtheorem{rem}[thm]{Remark}
 
 \numberwithin{equation}{section}

\newcommand{\NN}{\mathbb{N}}

\newcommand{\QQ}{\mathbb{Q}}
\newcommand{\RR}{\mathbb{R}}

\newcommand{\sgn}{\mathop{{sgn}}}

\begin{document}
\title[Characterization of Legendre curves in quasi-Sasakian pseudo-metric 3-manifolds]
{Characterization of Legendre curves in quasi-Sasakian pseudo-metric 3-manifolds}
\author[K. Srivastava, K. Sood, S. K. Srivastava] {K. Srivastava$^1$, K. Sood$^2$, S. K. Srivastava$^3$}
\address{$^{1,2,3}$Department of Mathematics,
                   Central University of Himachal Pradesh\br
                   Dharamshala-176215,
                   Himachal Pradesh,
	       INDIA}
\email{ksriddu22@gmail.com$^1$, soodkanika1212@gmail.com$^2$, \\ sachin@cuhimachal.ac.in$^3$.}

\thanks{K. Srivastava: supported by DST, Ministry of Science and Techonology, India through WOS-A vide their File No. SR/WOS-A/PM-20/2018. 
K. Sood: supported by DST, Ministry of Science and Techonology, India through JRF [IF160490] DST/INSPIRE/03/2015/005481} 
\begin{abstract}The main purpose of this paper is to present the spherical characterization of Legendre curves in $3$-dimensional quasi-Sasakian pseudo-metric manifolds. Furthermore, null Legendre curves are also characterized in this class of manifolds. 

\end{abstract}
\keywords{Contact pseudo-metric manifold, Frenet frame, Legendre curve, Spherical curve}
\subjclass{53A55, 53B25, 53C25, 53D15}
\maketitle
\section{Introduction}
The Legendre curves play a fundamental role in $3$-dimensional contact geometry. Let $(M;\varphi,\xi,\eta, g)$ be an almost contact metric $3$-manifold. Then an integral curve of the contact distribution $\ker\eta=\{X\in\Gamma(\mathcal{T}M)\,\vert\,\eta(X)=0\}$ is known as Legendre curve; $\Gamma\left(\mathcal{T}M\right)$ being the section of tangent bundle $\mathcal{T}M$ of $M$ \cite{BD,DEB}. The authors in \cite{camci, MB} have studied these curves in $3$-dimensional Sasakian pseudo-metric manifolds. In the present study, we will focus on some specific Legendre curves in quasi-Sasakian pseudo-metric $3$-manifolds. \\
From \cite{SB}, it is known that a regular curve $\upsilon:I\to\mathbb{E}^3$ in a $3$-dimensional Euclidean space $\mathbb{E}^3$ is said to be a spherical curve if and only if 
\begin{align}\label{de}
\frac{\tau(s)}{\kappa(s)}+\left(\frac{1}{\kappa(s)}\left(\frac{1}{\tau(s)}\right)'\right)'=0,
\end{align}
where curvature $\kappa(s)$ and torsion $\tau(s)$ of $\upsilon$ are non-zero smooth functions on open interval $I$, and the prime $'$ denotes differentiation with respect to arc length parameter $s$. The integral representation of Eq. \eqref{de} is given by 
\begin{align}\label{ie}
\kappa(s)^{-1}=C\cos\left(\int_{s_0}^{s}\tau(t)\,dt\right)+D\sin\left(\int_{s_0}^{s}\tau(t)\,dt\right).
\end{align}
In this study, we present spherical characterization of Legendre curves in quasi-Sasakian pseudo metric $3$-manifolds. Moreover, we obtain that null Legendre curves in this class of manifolds are geodesic.

\section{Preliminaries}\label{pre}
Let $M$ be a $C^\infty$ and paracompact $(2n+1)$-manifold, then $M$ is said to be an almost contact manifold if its structure group ${\rm GL}_{2n+1}\RR$ of tangent bundle $\mathcal{T}M$ is reducible to ${\rm U}(n)\times\{1\}$. This is equivalent to the existence of $(\varphi,\xi,\eta)$-structure satisfying
\begin{align}\label{2.1}
\varphi^2 = -\mathcal{I} +\eta\otimes\xi,\quad\eta\left(\xi\right)=1,
\end{align}
where $\varphi$ is an endomorphism, $\xi$ is a Reeb vector field, $\eta$ is a global $1$-form such that $\eta\wedge (d\eta)^n\ne 0$, $d$ being the exterior differential operator and $\mathcal{I}$ is the identity of $\mathcal{T}M$. We call $\eta$ the contact form. Eq. \eqref{2.1} yields $\varphi\xi = 0$, $\eta \circ\varphi =0$ and ${\rm rank}(\varphi)= 2n$. A pseudo-metric $g$ is said to be compatible with $(\varphi,\xi,\eta)$-structure if 
\begin{align}\label{2.3}
g\left(\varphi\cdot, \varphi\cdot \right)=g(\cdot, \cdot) -\varepsilon \eta(\cdot)\eta(\cdot),
\end{align}
where $\varepsilon^2=1$ and signature of $g$ is either $(2p, 2n-2p +1)$ or $(2p+1, 2n-2p)$ according as Reeb vector field is  timelike or spacelike respectively. Such $(M;\varphi,\xi,\eta, g)$ is called an almost contact pseudo-metric $(2n+1)$-manifold. Here, $\eta (X)= \varepsilon g(X, \xi)$ and $g(\xi,\xi)=\varepsilon$. Thus $\xi$ is never lightlike. The fundamental $2$-form $\Phi$ is defined by
$\Phi (\cdot, \cdot) = \varepsilon g(\cdot, \varphi \cdot)$.\\

On the direct product manifold $M\times\mathbb{R}$ : $\left(X, \ell\frac{d}{dt}\right)$ is an arbitrary tangent vector, where  $\ell$ is a smooth function on $M\times\mathbb{R}$, $t$ being standard coordinate on $\mathbb{R}$  and $X\in\Gamma(\mathcal{T}M)$. The almost complex structure $J$ on this direct product is defined by
\begin{align*}J\left(X, \ell\frac{d}{dt}\right)=\left(\varphi X - \ell\xi, \eta(X)\frac{d}{dt}\right).\end{align*} 
Then $M$ is said to be normal if and only if $J$ is integrable. Equivalently, $M$ is normal if and only if 
\begin{align}\label{2.7}
N_{\varphi}(X, Y) + 2 d\eta(X, Y)\xi = 0,
\end{align}
where $N_{\varphi}$ denotes the Nijenhuis torsion of endomorphism $\varphi$ which is given as follows:
$$N_{\varphi}(X, Y) = [\varphi X, \varphi Y] - \varphi[\phi X, Y]+ \varphi^2[X, Y]  - \varphi[X, \varphi Y],$$
for any tangent vectors $X$, $Y$ on $M$ (cf, \cite{DEB, DP2019, GC2010, KS}).\\
In this paper, we confine to dimension three. Following \cite{ZO}, we give certain results relevant to this case. Let $M$ be an almost contact pseudo-metric $3$-manifold, then we receive that
\begin{align}\label{nablax}
(\nabla_{X}\varphi)Y = - \eta(Y)\varphi\nabla_{X}\xi+\varepsilon g(\varphi\nabla_{X}\xi, Y)\xi ,\ X\in\Gamma (\mathcal{T}M).
\end{align}
\begin{prop}\label{equiv}
In an almost contact pseudo-metric $3$-manifold $M$, we have the following mutually equivalent conditions:
\begin{itemize}
\item[(i)] $M$ is normal,
\item[(ii)] $\nabla_{\varphi X}\xi=\varphi\nabla_X \xi$,
\item[(iii)] $\nabla_{X}{\xi} = - \varepsilon\alpha \varphi X + \varepsilon\beta \left(X-\eta\left(X\right)\xi \right)$,\end{itemize}
where $X\in\Gamma (\mathcal{T}M)$, $\alpha$ and $\beta$ being smooth functions on $M$ for which we have that
\begin{align}\label{3.14}
2 \alpha = {\rm trace}\left\lbrace X \rightarrow \varphi\nabla_{X}\xi\right\rbrace, \,\,\, 2\beta = {\rm trace}\left\lbrace X \rightarrow\nabla_{X}\xi\right\rbrace.
\end{align}
\end{prop}
\noindent Using Eq.\eqref{nablax} and above proposition, we have 
\begin{align}\label{3.1}
(\nabla_{X}\varphi)Y = \beta(g(\varphi X, Y)\xi - \varepsilon\eta(Y)\varphi X)+\alpha(g(X, Y)\xi - \varepsilon\eta (Y)X).
\end{align}
Let us denote by $\mathbb{N}^3$ a normal almost contact pseudo-metric $3$-manifold. Analogous to \cite{Inoguchi15}, $\mathbb{N}^3$ is called a quasi-Sasakian pseudo-metric $3$-manifold (denoted by $\mathbb{Q}^3$) if $\beta=0$ and $\xi(\alpha)=0$, in particular Sasakian if $\beta=0$ and $\alpha=1$. Below, we give an example of $\mathbb{N}^3$.
\begin{ex}
Let $(x,y,z)$ be the standard Cartesian coordinates on $\RR^3$ and consider the $1$-form $\eta=2ydx+dz$. We put $\xi=\partial_3$ and consider the endomorphism $\varphi$ defined by 
$\varphi\partial_1=\partial_2,\ \varphi\partial_2=2y\partial_3-\partial_1,\ \varphi\partial_3=0$, where $\partial_1=\frac{\partial}{\partial x}$, $\partial_2=\frac{\partial}{\partial y}$ and $\partial_3=\frac{\partial}{\partial z}$. Then it follows that $\eta(\xi)=1$ and $\varphi^2=-\mathcal{I}+\eta\otimes\xi$. Hence $(\varphi,\xi,\eta)$ is an almost contact structure on $\RR^3$. By straightforward computations,we obtain that 
\begin{align*}
N_{\varphi}(X, Y) \left(\partial_i ,\partial_j\right)+2d\eta\left(\partial_i ,\partial_j\right)\xi=0,\quad i,j\in\{1,2,3\}
\end{align*}
which implies that the structure is normal.
 
\noindent Let $\mathcal{N}^3_{\varepsilon}:=\RR^3$ and consider a normal almost contact structure $(\varphi,\xi,\eta)$ restricted to $\mathcal{N}^3_{\varepsilon}$. Next, we consider the metric tensor $g=\exp(2z)(dx^2+dy^2)+\varepsilon\eta\otimes\eta$, with matrix representation with respect to the standard basis
 \begin{align*}
 g[(\partial_i, \partial_j)]=\begin{bmatrix}4\varepsilon y^2 +\exp(2z)&0& 2\varepsilon y \\ 0&\exp(2z)&0\\ 2\varepsilon y&0&\varepsilon \end{bmatrix},
 \end{align*}
where $\varepsilon^2=1$. Using the endomorphism $\varphi$ and the metric $g$, we obtain $g\left(\varphi X, \varphi Y\right)=g(X, Y) -\varepsilon \eta(X)\eta(Y)$ and $\eta(X)=\varepsilon g(X,\xi)$, and hence that $(\mathcal{N}^3_{\varepsilon};\varphi,\xi,\eta,g)$ is a normal almost contact pseudo-metric $3$-manifold. For the Levi-Civita connection $\nabla$ with respect to this metric, we have
\begin{align*}
\nabla_{\partial_1}\partial_1=&2y\partial_1 -4\varepsilon y\exp(-2z)\partial_2 -\varepsilon(4\varepsilon y^2+\exp(2z))\partial_3,\\ 
\nabla_{\partial_1}\partial_2 =& 2\varepsilon y\exp(-2z)\partial_1 +(-4\varepsilon y^2+\exp(2z))\exp(-2z)\partial_3=\nabla_{\partial_2}\partial_1, \\ 
\nabla_{\partial_1}\partial_3=& \partial_{1}-\varepsilon \exp(-2z)\partial_2 -2y\partial_3=\nabla_{\partial_3}\partial_1,\\ 
\nabla_{\partial_2}\partial_2=&2y\partial_1 -\varepsilon(4\varepsilon y^2+\exp(2z))\partial_3, \\
 \nabla_{\partial_2}\partial_3=&\varepsilon \exp(-2z)\partial_1+\partial_2 -2\varepsilon y\exp(-2z)\partial_3=\nabla_{\partial_3}\partial_2, \nabla_{\partial_3}\partial_3=0.\nonumber
\end{align*}
Using Eq. \eqref{3.1} and the above expressions, we find $\alpha=\exp(-2z)$ and $\beta=\varepsilon$. \end{ex}


\section{Legendre curves}\label{lc1}
Let $ M$ be an almost contact pseudo-metric $3$-manifold with Levi-Civita connection $\nabla$ and $\upsilon : I\longrightarrow M$ a unit speed curve parametrized by the arc-length $s$ in $M$,  $I$ being an open interval, then the Frenet frame (or Frenet $3$-frame) $\{T\coloneqq\upsilon',\ N,\ B\}$ of $\upsilon$ satisfies the following {\it Frenet-Serret} formulas:  
\begin{align}
  \nabla_{T}T = \kappa N, \,\,\, \nabla_{T}N = -\kappa T + \varepsilon \tau B \,\,\, {\rm{and}} \,\,\, \nabla_{T}B = -\tau N,\label{4.1}
\end{align}
where $\kappa=\vert\nabla_{T}{T}\vert$  and $\tau$ are the geodesic curvature and geodesic torsion of $\upsilon$, respectively. The vector fields $T$, $B$ and $N$ are known as the tangent, binormal and principal normal vector fields of $\upsilon$, respectively.
If $\nabla_{\upsilon '}\upsilon ' = 0$ then $\upsilon$ is said to be a geodesic and non-geodesic if $\kappa>0$ everywhere on $I$.  \\
\noindent If $\upsilon$ is an integral curve of the contact distribution $\mathcal{D}={\rm Ker}\, \eta$, equivalently, $\eta(\upsilon ') = 0$ then we say that $\upsilon$ is a {\it Legendre curve} in $M$. These curves have been intensively studied by several authors (see \cite{BD, camci, CM14, CM15, Inoguchi, JW07, MB}).\\

\subsection{Legendre curves in $\NN^3$} 
Let $\upsilon:I\rightarrow \NN^3$ be a Legendre curve in $\NN^3$. Then $\{\upsilon',\varphi\upsilon',\xi\}$ is an orthonormal frame along $\upsilon$, where $g(\upsilon',\upsilon')=g(\varphi\upsilon',\varphi\upsilon')=1$ and $g(\xi,\xi)=\varepsilon$. Let us put $\nabla_{T}T=a_1\upsilon'+\vartheta\varphi\upsilon'+c_1\xi,$ where $a_1,\vartheta$ and $c_1$ are some smooth functions on $I$. Then we find that $a_1=g(\nabla_{\upsilon'}\upsilon', \upsilon')=0$, $\vartheta=g(\nabla_{\upsilon'}\upsilon', \varphi\upsilon')$ and $-\varepsilon c_1=g(\upsilon',\nabla_{\upsilon'}\xi)$. Employing Eqs. \eqref{2.3} and \eqref{3.1}, we have $g(\upsilon',\nabla_{\upsilon'}\xi)=\varepsilon\beta$. Thus, 
\begin{align}\label{nablaTT}
\nabla_{T}T=-\beta \xi+\vartheta \varphi \upsilon',
\end{align}
which yields $\kappa=\sqrt{\vert\vartheta^2+\varepsilon\beta^2\vert}$. If $\kappa>0$ everywhere on $I$ then differentating $N=\frac{1}{\kappa}\nabla_T T=-\frac{\beta}{\kappa}\xi+\frac{\vartheta}{\kappa}\varphi \upsilon'$ along $\upsilon$ and using Eqs. \eqref{3.1} and \eqref{nablaTT}, we obtain that  
 \begin{align}
 \nabla_{T}N=\frac{c}{\kappa}\left(\vartheta \xi+\epsilon \beta \varphi \upsilon'\right)-\kappa\upsilon', \label{nablaN}
 \end{align}
where $c=\alpha+\frac{\beta \vartheta'-\beta' \vartheta}{\epsilon\beta^2+\vartheta^2}$. In light of \eqref{4.1},  we have from Eq. \eqref{nablaN} that $\varepsilon\tau B=\frac{c}{\kappa}(\vartheta \xi+\epsilon \beta \varphi \gamma')$, which provides $\tau=\vert c\vert$. This leads to the following result which generalize Theorem 1 of \cite{JW07} to the case of spacelike or timelike Reeb vector field.
\begin{prop}\label{prop}
Let $\upsilon : I\longrightarrow \NN^3$ be a non-geodesic Legendre curve in $\NN^3$. Then its curvature $\left(\kappa\right)$ and torsion $\left(\tau\right)$ are given by
\begin{align}\label{4.2}
\kappa=&\sqrt{\vert\vartheta^2+\varepsilon\beta^2\vert},\\
\label{4.3} \tau=&\left\arrowvert\alpha +\left(\frac{\beta\vartheta'-\beta'\vartheta}{\kappa^2}\right)\right\arrowvert,
\end{align}
where $\vartheta$ is a smooth function defined by $\vartheta=g\left(\nabla_{\gamma'}\gamma', \varphi\gamma'\right)$
on $I$, and $\alpha$, $\beta$ being the same as in \eqref{3.14}.
\end{prop}
\begin{rem}
\begin{itemize} 
\item[(a)]In $\NN^3$, the decomposition of Reeb vector field $\xi$ in the Frenet frame of a non-geodesic Legendre curve is given by
\begin{align}
\xi=\frac{\epsilon}{\sqrt{\vert\vartheta^2+\varepsilon\beta^2\vert}}\left\{-\beta N+\vartheta B\right\}.\label{xi}
\end{align}
\item[(b)] Let $\upsilon$ be a non-geodesic Legendre curve in $\NN^3$ with $\xi$ is parallel to binormal $B$. Then $\NN^3$ is a quasi-Sasakian pseudo-metric $3$-manifold.
\item[(c)] Let $\upsilon$ be a non-geodesic Legendre curve in $\NN^3$ such that $\xi$ is parallel to principal normal $N$ then it has torsion $\vert\alpha\vert$. In particular, the Legendre curves of a $\beta$-Kenmotsu manifold with $\xi$ is spacelike and parallel to $N$ are circles \cite{CMM}. 
\end{itemize}
\end{rem} 
\begin{ex}
\noindent If we define a curve $\upsilon:I\subset\RR\rightarrow\mathcal{N}^3_{\varepsilon}$ by $\upsilon(s)=\left(\upsilon_1(s),\upsilon_2(s),\upsilon_3(s)\right)$. Then $\upsilon$ is a Legendre curve if and only if
\begin{align}\label{con}
\begin{cases}
2\upsilon_1'\upsilon_2+\upsilon_3'=0,\\
\upsilon_1'^2+\upsilon_2'^2=\exp(-2\upsilon_3).\end{cases} 
\end{align}
It can be easily seen that
\begin{align*}
\upsilon'=\upsilon_1'\partial_1+\upsilon_2'\partial_2-2\upsilon_1'\upsilon_2\partial_3\ {\rm and}\
\varphi\upsilon'=-\upsilon_2'\partial_1+\upsilon_1'\partial_2+2\upsilon_2\upsilon_2'\partial_3.
\end{align*}
Below, we present certain Legendre curves in $\mathcal{N}^3_{\varepsilon}$.
\begin{itemize}
\item[(i)] Let us consider a curve 
\begin{align*}
\upsilon^1(s)=(1,s,0)
\end{align*}
in $\mathcal{N}^3_{\varepsilon}$, where $s\in I\subseteq\RR$ is arc length parameter of $\upsilon^1$. Then this curve satisfies Eq. \eqref{con}, so it is a Legendre curve in $\mathcal{N}^3$. For such a curve $\alpha(\upsilon(s))=1$, $\beta(\upsilon(s))=\varepsilon$, $\vartheta=g\left(\nabla_{\upsilon'}\upsilon',\varphi\upsilon'\right)=-2s$, $$\kappa=\sqrt{\vert1+4\varepsilon s^2\vert} \ {\rm and}\ \tau=\left\vert1-\frac{2\varepsilon}{\kappa^2}\right\vert.$$
\item[(ii)] Let us consider a curve 
\begin{align*}
\upsilon^2(s)=(-\ln s,1/2,\ln s)
\end{align*}
in $\mathcal{N}^3_{\varepsilon}$, where $s>0$ is arc length parameter of $\upsilon^2$. Then this curve satisfies Eq. \eqref{con}, so it is a Legendre curve in $\mathcal{N}^3$. For such a curve $\alpha(\upsilon(s))=1/{s^2}$, $\beta(\upsilon(s))=\varepsilon$, $\vartheta=0$, $$\kappa=1 \ {\rm and}\ \tau=\frac{1}{s^2}.$$
\end{itemize}
Further, the spacelike and timelike situations of these Legendre curves are interesting. The Euclidean pictures of $\upsilon^1$ and $\upsilon^2$ are rendered in Fig. $1$ and Fig. $2$ respectively.
\begin{figure}[ht]
\begin{minipage}{0.3\textwidth}
\includegraphics[height=5 cm]{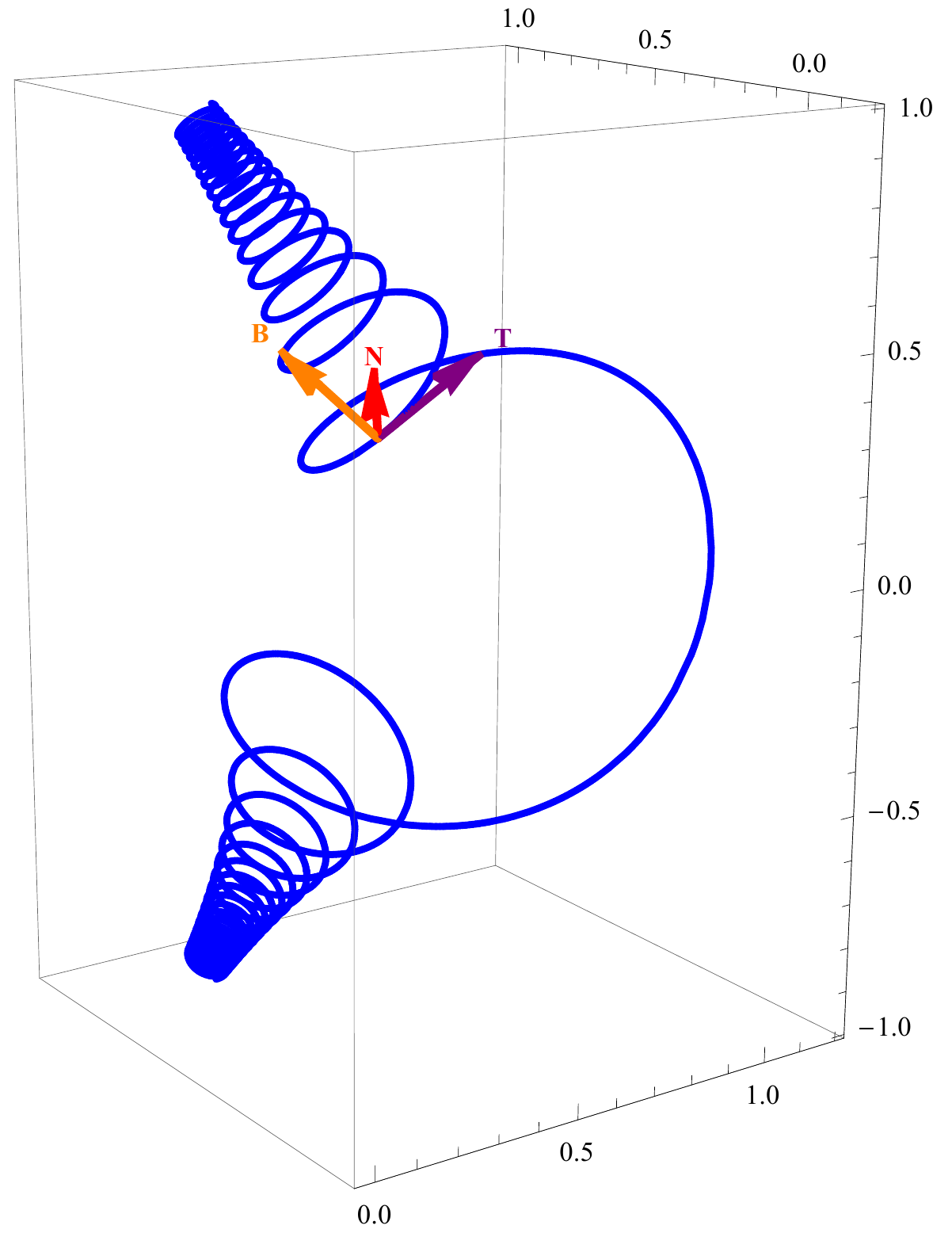}\\
\center(a) $\upsilon^1$ in $\mathcal{N}^3_{1}$ 
\end{minipage}
\hspace*{2 cm}
\begin{minipage}{0.3\textwidth}
\includegraphics[height=5 cm]{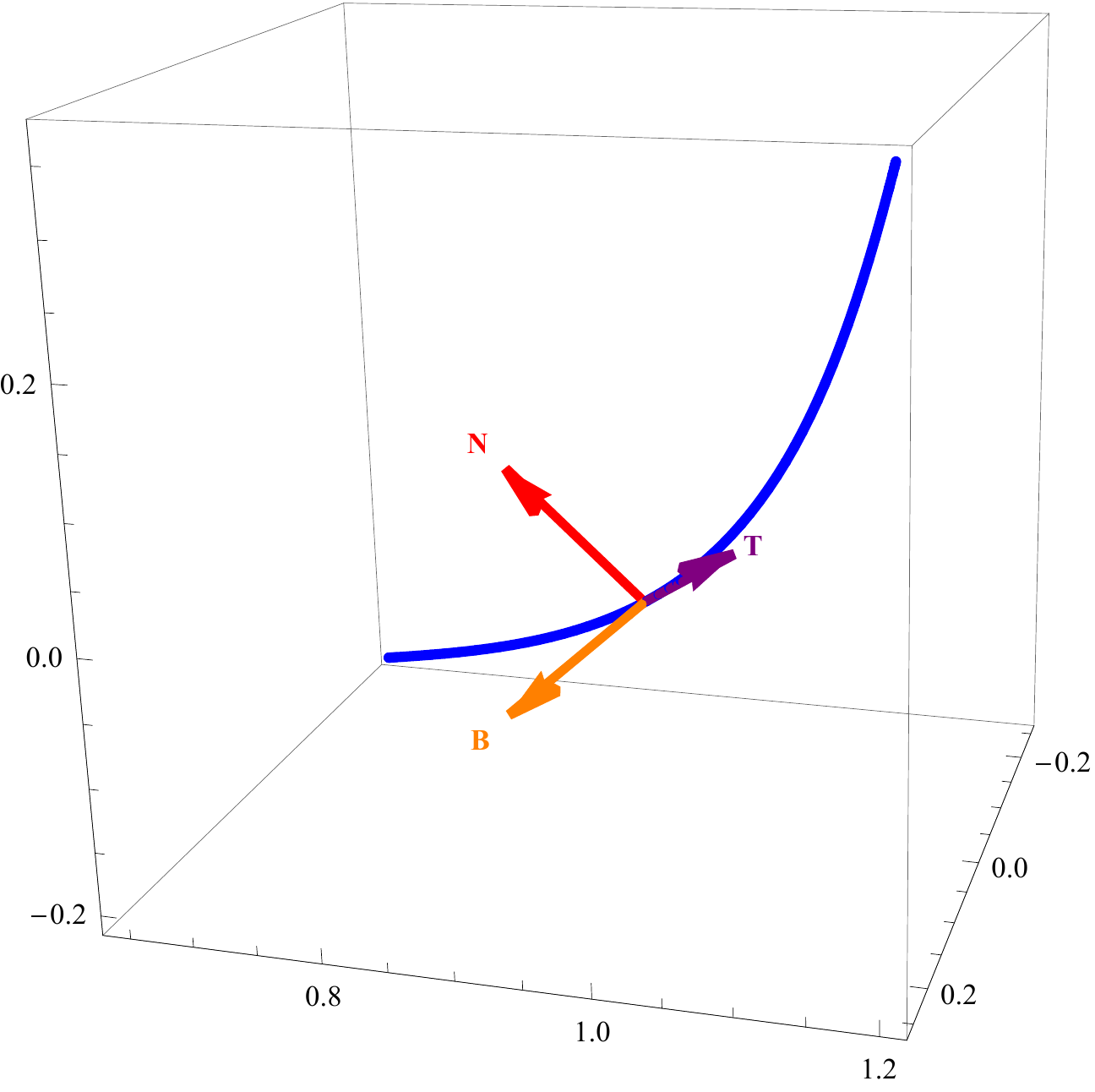}\\
\center(b) $\upsilon^1$ in $\mathcal{N}^3_{-1}$  
\end{minipage}
\caption{Legendre curve $\upsilon^1$ in $\mathcal{N}^3_{\varepsilon}$}
\end{figure}

\begin{figure}[ht]
\begin{minipage}{0.3\textwidth}
\includegraphics[height=5 cm]{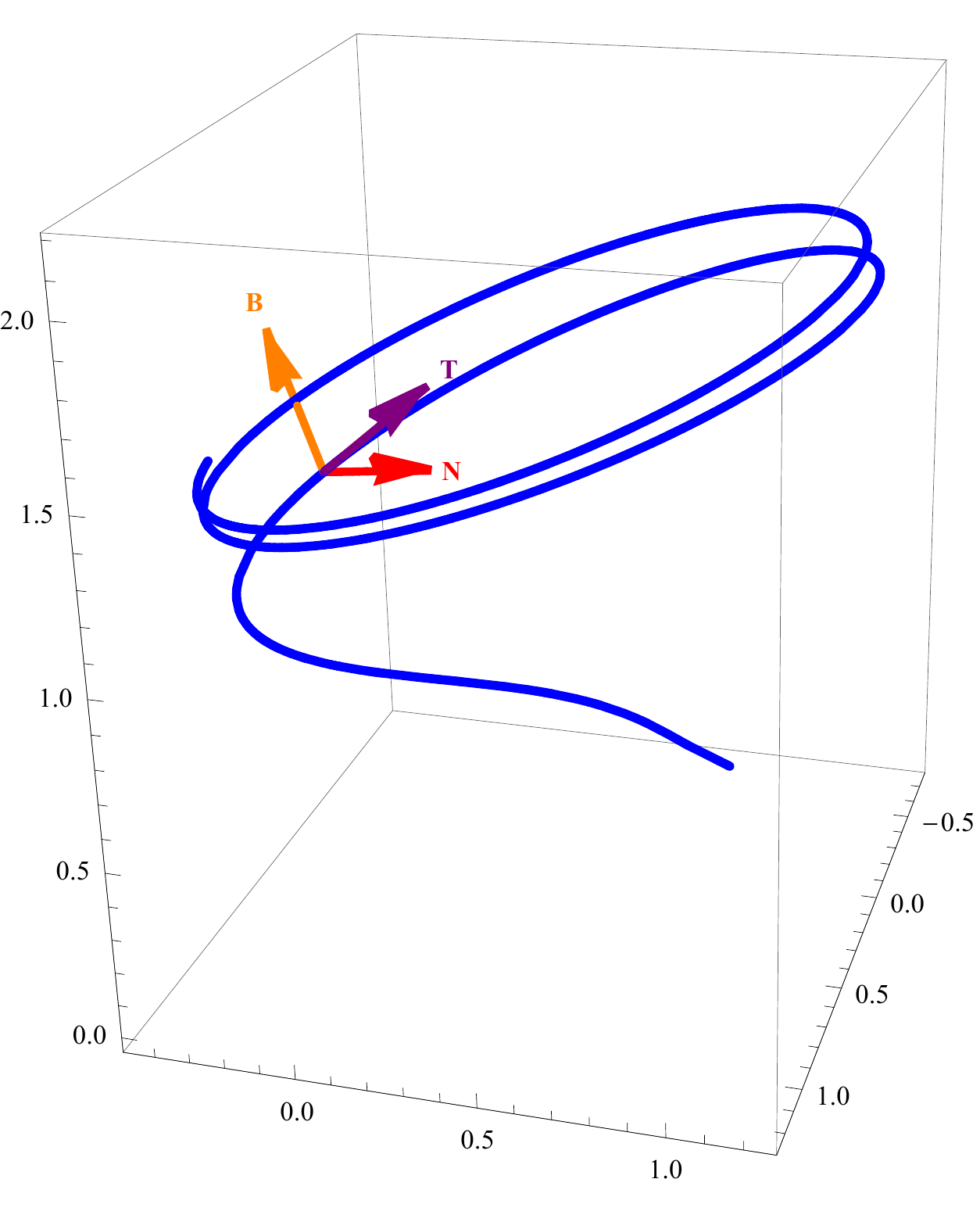}\\
\center(c) $\upsilon^2$ in $\mathcal{N}^3_{1}$ 
\end{minipage}
\hspace*{2 cm}
\begin{minipage}{0.3\textwidth}
\includegraphics[height=5 cm]{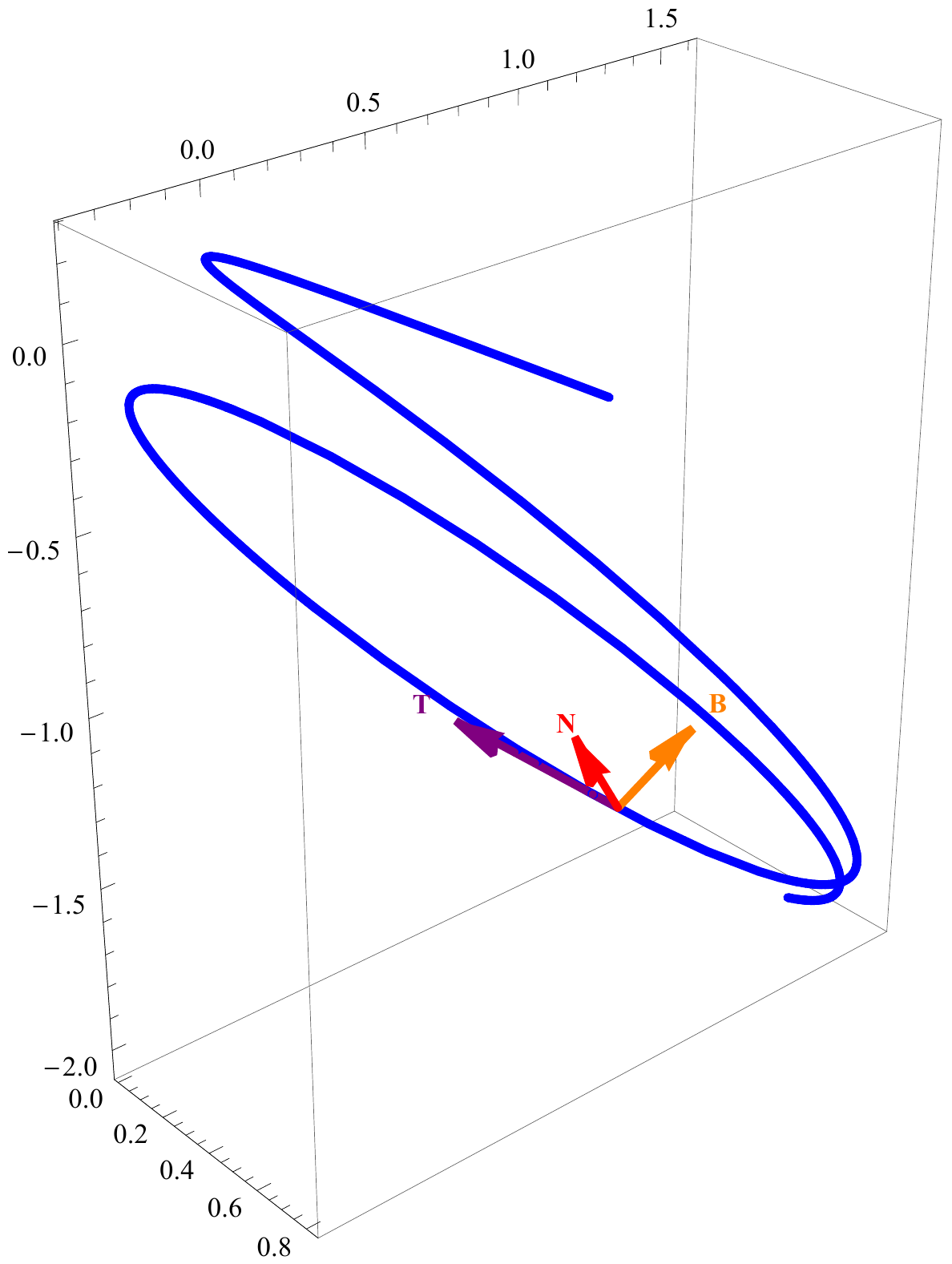}\\
\center(d) $\upsilon^2$ in $\mathcal{N}^3_{-1}$ 
\end{minipage}
\caption{Legendre curve $\upsilon^2$ in $\mathcal{N}^3_{\varepsilon}$}
\end{figure}
\end{ex}

\noindent Let $\upsilon:I\rightarrow \NN^3$ be a unit speed curve such that $\eta(\upsilon')=m$, where $m$ is a smooth function on $I$. Let us consider the vector fields $\upsilon', \varphi\upsilon'$ and $\xi$ for which, we have that $g(\upsilon', \upsilon')=1$, $g(\upsilon', \xi)=\varepsilon m$, $g(\varphi\upsilon', \varphi\upsilon')=1-\varepsilon m^2$, $g(\xi, \xi)=\varepsilon$ and $g(\upsilon',\varphi\upsilon')=g(\xi,\varphi\upsilon')=0$. Then, $\upsilon', \varphi\upsilon', \xi$ are linearly independent. Further, $\{\upsilon', \varphi\upsilon', \xi\}$ forms a basis of $\mathcal{T}_{\upsilon(s)}M^3$ for any $s\in I$ if and only if $1-\varepsilon m^2\ne 0.$ In this case, one can define orthonormal frame $\{V_1,V_2,V_3\}$ along $\upsilon$ as follows:
\begin{align}\label{frame}
V_1=\upsilon',\quad\quad V_2=\frac{\varphi\upsilon'}{\delta},\quad\quad V_3=\frac{\xi-\varepsilon m\upsilon'}{\delta} 
\end{align}
where $g(V_1, V_1)=g(V_2, V_2)=1$, $g(V_3, V_3)=\varepsilon$ and $\delta=\sqrt{\vert1-\varepsilon m^2\vert}$.

\noindent It is important to note that $\{\upsilon', \varphi\upsilon', \xi\}$ is linearly dependent if and only if either $\upsilon'=m\xi$ or $\upsilon'=m\xi+\varepsilon\varphi\upsilon'$. In reality, if these vector fields are linearly dependent then $\xi$ is spacelike and $\vert m\vert=1$. In this case, the curve $\upsilon$ is necessarily a geodesic. Therefore, for the non-geodesic curve $\upsilon$, we must have that neither $\varepsilon= 1$ nor $\vert m\vert=1$, that is $\delta\ne 0$. The decomposition of characteristic vector field with respect to the frame $\left\{V_1,V_2,V_3\right\}$ is given by 
\begin{align}
\xi=\varepsilon\left(mV_1+\delta V_3\right).\label{dxi}
\end{align}
 Below we present some preparatory lemmas for later use.
\begin{lem}
Let $\upsilon:I\rightarrow \NN^3$ be a non-geodesic unit speed curve in $\NN^3$ such that $\eta(\upsilon')=m$. Then the covariant derivative of $V_1, V_2, V_3$ along $\upsilon$ are given by:
\begin{align}
\nabla_{\upsilon'}V_1=&\delta\vartheta_1 V_2-\left(\beta\delta-\frac{m'}{\delta}\right)V_3,\label{E1}\\
\nabla_{\upsilon'}V_2=&-\delta\vartheta_1V_1+(\alpha+m\vartheta_1)V_3,\label{E2}\\
\nabla_{\upsilon'}V_3=&\varepsilon\left( \left(\beta\delta-\frac{m'}{\delta}\right) V_1-(\alpha+m\vartheta_1)V_2\right),\label{E3}
\end{align}
where 
\begin{align}\label{abcm}
\delta=\sqrt{\vert1-\varepsilon m^2\vert} \ {\rm and}\ \vartheta_1=\frac{1}{\delta^2}g(\nabla_{\upsilon'}\upsilon',\varphi\upsilon').
\end{align}
\end{lem}
\begin{proof}
Let us express $\nabla_{\upsilon'}V_1=a_2V_1+b_2V_2+c_2V_3$, where $a_2,\ b_2, \ c_2$ are some smooth functions on $I$. Then $a_2=g(\nabla_{\upsilon'}V_1, V_1)=0$. 

\noindent Using Eq. \eqref{3.1}, we find 
\begin{align*}
-b_2&=g(V_1, \nabla_{\upsilon'}V_2)=\frac{1}{\delta}g((V_1,(\nabla_{\upsilon'}\varphi)\upsilon'+\varphi\nabla_{\upsilon'}\upsilon')\\
&=-\frac{1}{\delta}g(\varphi V_1, \nabla_{\upsilon'}\upsilon'))=-\delta\vartheta_1.
\end{align*}

\noindent Employing Proposition \eqref{equiv}, we have 
\begin{align*}
-\varepsilon c_2&=g(V_1,\nabla_{\upsilon'}V_3)=\frac{1}{\delta}g(V_1,\nabla_{\upsilon'}\xi-\varepsilon m'\upsilon')\\
&=\frac{1}{\delta}g(V_1,\varepsilon\beta(\upsilon'-m\xi)-\varepsilon m'\upsilon')
=\frac{\varepsilon}{\delta}(\beta \delta^2-m'). 
\end{align*}
In light of above expressions, we obtain that $a_2=0$, $b_2=\delta\vartheta_1$, $c_2=-\frac{\varepsilon}{\delta}(\beta \delta^2-m')$ and this leads to Eq. \eqref{E1}. Similarly, we compute $\nabla_{\upsilon'}V_2$ and $\nabla_{\upsilon'}V_3$. 
\end{proof}

\begin{lem}
Let $\upsilon:I\rightarrow \NN^3$ be a non-geodesic unit speed curve in $\NN^3$ such that $\eta(\upsilon')=m$. Then its curvature and torsion are given by
 \begin{align}\label{kt1}
  \begin{cases}
 \kappa=&\delta\sqrt{\left\vert\vartheta_1^2+\varepsilon\left(\beta-\frac{m'}{\delta^2}\right)^2\right\vert},\\ 
\tau=&\left\vert \alpha+m\vartheta_1+\frac{\left(\beta\vartheta_1'-\beta'\vartheta_1\right)-2\left(\frac{m'\vartheta_1'}{\delta^2}\right)+\left(\frac{m'\vartheta_1}{\delta^2}\right)'}{\vartheta_1^2+\varepsilon\left(\beta -\frac{m'}{\delta^2}\right)^2}\right\vert,
 \end{cases}
 \end{align}
 where $\delta$ and $\vartheta_1$ are the same as in Eq. \eqref{abcm}.
\end{lem}
\begin{proof}
\noindent Using Eq. \eqref{E1} and computing the length of $\nabla_{\upsilon'}\upsilon'$, we receive the curvature function $\kappa$. By the virtue of Eqs. \eqref{4.1} and \eqref{E1}, we find 
\begin{align}\label{N}
N=\frac{\delta\vartheta_1}{\kappa}V_2-\left(\frac{\beta\delta}{\kappa}-\frac{m'}{\kappa\delta}\right)V_3.
\end{align}

\noindent Let us denote $a_2=\delta\vartheta_1$ and $b_2=\left(\beta\delta-\frac{m'}{\delta}\right)$. Now, we compute
\begin{align}
\nabla_{\upsilon'}N=&\left(\frac{a_2}{\kappa}\right)'V_2+\left(\frac{a_2}{\kappa}\right)\nabla_{\upsilon'}V_2-\left(\frac{b_2}{\kappa}\right)'V_3-\left(\frac{b_2}{\kappa}\right)\nabla_{\upsilon'}V_3\nonumber\\
=&\left(-\frac{{a_2}^2}{\kappa}-\frac{\varepsilon {b_2}^2}{\kappa}\right)V_1+\left(\frac{{a_2}'\kappa-a_2{\kappa}'}{\kappa^2}+\frac{\varepsilon b_2}{\kappa}\left(\alpha+\frac{a_2 m}{\delta}\right)\right)V_2 \nonumber\\
& +\left(-\frac{{b_2}'\kappa-b_2{\kappa}'}{\kappa^2}+\frac{a_2}{\kappa}\left(\alpha+\frac{a_2 m}{\delta}\right)\right)V_3\nonumber\\
=&-\kappa V_1+\frac{\varepsilon b_2}{\kappa}\left(\alpha+\frac{a_2 m}{\delta}+\frac{{a_2}'b_2-a_2{b_2}'}{\kappa^2}\right)V_2+\frac{a_2}{\kappa}\left(\alpha+\frac{a_2 m}{\delta}+\frac{{a_2}'b_2-a_2{b_2}'}{\kappa^2}\right)V_3\nonumber\\
=&-\kappa V_1+\left(\alpha+\frac{a_2 m}{\delta}+\frac{{a_2}'b_2-a_2{b_2}'}{\kappa^2}\right)\left(\frac{\varepsilon b_2V_2+a_2V_3}{\kappa}\right).\label{t}
\end{align}
In light of the Frenet-Serret formulas, we have from Eq. \eqref{t} that $$\tau=\left\vert \alpha+\frac{a_2 m}{\delta}+\frac{{a_2}'b_2-a_2{b_2}'}{\kappa^2}\right\vert.$$  Using the values of $a_2$ and $b_2$ in the above expression, we find the torsion function of $\upsilon$ as mentioned in Eq. \eqref{kt1}. This completes the proof.
\end{proof}
\begin{thm}
Let $\upsilon:I\rightarrow \NN^3$ be a non-geodesic unit speed curve in $\NN^3$ such that $\eta(\upsilon')=m$. Then, the decomposition of Reeb vector field $\xi$ in the Frenet frame $\{T:=\upsilon',N,B\}$ of $\upsilon$ is given by
\begin{align}\label{decomp}
\xi={\varepsilon}(m T+\eta(N)N)+\eta(B)B,
\end{align}
 where $\eta(N)=(m'-\beta\delta^2)/\kappa$, $\eta(B)=(\varepsilon\sgn(\tau)\delta^2\vartheta_1)/\kappa$, $\delta$ and $\vartheta_1$ are the same as in Eq. \eqref{abcm}.
\end{thm}
\begin{proof}\ Let us put $\xi=pT+qN+rB$, where $p,q$ and $r$ are smooth functions on $I$. Then $p=g(\xi,T)=\varepsilon m$, $q=g(\xi,N)=\varepsilon\eta(N)$ and $r=\varepsilon g(\xi,B)=\eta(B)$. The values of $\eta(N)$ and $\eta(B)$ directly follows from Eqs. \eqref{4.1}, \eqref{dxi}, \eqref{N} and \eqref{t}. This completes the proof.

\end{proof}
\noindent In view of Eq. \eqref{decomp}, we find that $\eta(B)^2+\varepsilon\eta(N)^2=\delta^2$. This leads to the following remark.
\begin{rem}
For a non-geodesic unit speed curve $\upsilon$ in $\NN^3$ such that $\eta(\upsilon')=m$, we have
\begin{align*}
\delta^2-\varepsilon\eta(N)^2\ge 0.
\end{align*}
\end{rem}

\subsection{Legendre curves in $\QQ^3$} 
Let $\upsilon$ be a non-geodesic Legendre curve in a quasi-Sasakian pseudo-metric $3$-manifold $\QQ^3$. Then by the virtue of proposition \ref{prop}, curvature and torsion of $\upsilon$ are given by
\begin{align}\label{quasikt}
\kappa=\vartheta,\ \tau=\vert\alpha\vert,
\end{align}
where $\vartheta=g\left(\nabla_{\upsilon'}\upsilon', \varphi\upsilon'\right)>0$ and $\alpha$ being the same as in \eqref{3.14}.
\begin{thm} \label{thm Br}
For a non-geodesic unit speed curve $\upsilon$ in $\QQ^3$ with $\eta(\upsilon')=m$. If $ \tau=\vert\alpha\vert$ and at one point $m=m'=m'' =0$, where $\alpha$ is non-zero at every point of $M$. Then $\upsilon$ is a Legendre curve.
\end{thm}
\begin{proof}
For spacelike vector $\xi$, this is Theorem 3 of \cite{JW07}, so we concentrate on the timelike situation. 

\noindent In this class of manifold, we have from Eq. \eqref{kt1} that
\begin{align}\label{tauquasi}
\tau=&\left\vert \alpha+m\vartheta_1+\left(\left(\frac{m'\vartheta_1}{\delta^2}\right)'-2\left(\frac{m'\vartheta_1'}{\delta^2}\right)\right)\left(\vartheta_1^2-\frac{m'^2}{\delta^4}\right)^{-1}\right\vert,
\end{align}
where $\delta=\sqrt{1+ m^2}\ge 1$ and $\vartheta_1\delta^2=g(\nabla_{\upsilon'}\upsilon',\varphi\upsilon')$. Since $m(s_0)=m'(s_0)=m''(s_0)=0$ at certain $s_0\in I$, therefore $\tau(s_0)=\vert\alpha(s_0)\vert\ne 0$. Further, $\vartheta_1(s_0)\ne 0$ since $\upsilon$ is non-geodesic. Taking these into consideration and using $\tau=\vert\alpha\vert$, Eq. \eqref{tauquasi} yields
\begin{align}\label{vartheta1m}
\left(m''\vartheta_1-m'\vartheta_1'\right)-3mm'\delta^{-2}\vartheta_1+m\delta^2\vartheta_1^3=0.
\end{align}
Obviously $m=0$ is a solution and $m'=0$ implies that $m=0$. Hence we now assume that $m$ is not constant. Setting $t=m'/\vartheta_1$, Eq. \eqref{vartheta1m} becomes
$\frac{d}{dm}\left(t^2\right)-\frac{6m}{\delta^2}t^2+2m\delta^2=0$.
Integration gives ${m'}^2=(1/3+C\delta^3)\delta^3\vartheta_1^2$, where $C$ is constant of integration. Using initial conditions, we have $C=-1/3$. Thus ${m'}^2=-(\delta^3-1)\delta^3\vartheta_1^2/3$. Since $\delta\ge 1$ and $\vartheta_1\ne 0$, we have $m=0$, a contradiction. This completes the proof.
\end{proof}

\begin{ex}
Let $\mathcal{Q}^3_{\alpha}:=\RR_{+}\times\RR^2$ with standard Cartesian coordinates $(x,y,z)$, $\eta=-2xdy+dz$, $\xi=\partial_3$, $\varphi$ defined by: $\varphi\partial_1=\partial_2+2x\partial_{3},\ \varphi\partial_2=-\partial_1,\ \varphi\partial_3=0$ and metric tensor $g=x^{2}(dx^2+dy^2)+\varepsilon\eta\otimes\eta$, where $\partial_1=\frac{\partial}{\partial x}$, $\partial_2=\frac{\partial}{\partial y}$ and $\partial_3=\frac{\partial}{\partial z}$. Then $(\mathcal{Q}^3_{\alpha};\varphi,\xi,\eta,g)$  is a normal almost contact pseudo-metric $3$-manifold.  For the Levi-Civita connection $\nabla$ with respect to this metric, we find
\begin{align}\label{levi}
&\nabla_{\partial_1}\partial_1=\frac{1}{x}\partial_1, \nabla_{\partial_2}\partial_2=-\frac{(1+4\varepsilon)}{x}\partial_1, \nabla_{\partial_3}\partial_3=0, \nabla_{\partial_2}\partial_3=\nabla_{\partial_3}\partial_2=\frac{\varepsilon}{x^{2}}\partial_1, \nonumber
\\
&\nabla_{\partial_1}\partial_2=\nabla_{\partial_2}\partial_1=\frac{(2\varepsilon+1)}{x}\partial_2+(1+4\varepsilon)\partial_3, \nabla_{\partial_1}\partial_3=\nabla_{\partial_3}\partial_1=-\frac{\varepsilon}{x^{2}}\partial_2-\frac{2\varepsilon}{x}\partial_3.
\end{align}
In view of Eqs. \eqref{3.1} and \eqref{levi}, we have $\alpha=1/x^{2}$ and $\beta=0$. Thus $\mathcal{Q}^3_{\alpha}$ is a quasi-Sasakian $3$-manifold.
\end{ex}
\noindent Using definition of Legendre curve, we find that $\upsilon:I\rightarrow\mathcal{Q}^3_{\alpha}$ defined by $\upsilon(s)=\left(\upsilon_1(s),\upsilon_2(s),\upsilon_3(s)\right)$ is a Legendre curve if and only if
\begin{align}\label{k1}
\begin{cases}
\upsilon_3'=2\upsilon_1\upsilon'_2,\\
\upsilon_1'^2+\upsilon_2'^2=\upsilon_1^{-2}.\end{cases} 
\end{align}
In this case, we have
\begin{align*}
\upsilon'=\upsilon_1'\partial_1+\upsilon_2'\partial_2+2\upsilon_1\upsilon'_2\partial_3\ {\rm and}\
\varphi\upsilon'=-\upsilon_2'\partial_1+\upsilon_1'\partial_2+2\upsilon_1\upsilon'_1\partial_3.
\end{align*}
From Eq. \eqref{k1}, we have  
\begin{align*}
\upsilon_1'=\upsilon_1^{-1}\cos\phi(s), \upsilon_2'=\upsilon_1^{-1}\sin\phi(s), \upsilon_3'=2\sin\phi(s),
\end{align*}
where $\phi\in C^{\infty}(I)$. Thus, we have the general expression of Legendre curves in $\mathcal{Q}^3_{\alpha}$ as follows:
\begin{prop}\label{lq3}
Let $\upsilon:I\rightarrow\mathcal{Q}^3_{\alpha}$ be a non-geodesic Legendre curve, then $\upsilon$ is given by
\begin{align}\label{k2}
\upsilon(s)=\left(\int_{s_0}^{s}{\mu^{-1}}(t)\zeta(t)dt, 2\int_{s_0}^{s}\sin\phi(s) \right),
\end{align}
where $\phi$ being smooth function on $I$, $\zeta(s)=(\cos\phi(s), \sin\phi(s))$ is any parametrization of circle $\mathbb{S}^1$ and $\mu$ is a non-zero smooth function defined by 
$\mu^2(s)=2\int_{s_0}^{s}\cos\phi(t)dt$.
\end{prop}
\begin{thm}
The curvature and torsion of non-geodesic Legendre curve \eqref{k2} in $\mathcal{Q}^3_\alpha$ are given by
\begin{align}\label{kt3}
\kappa=\phi'+\frac{\sin\phi}{\mu^2},
\tau=\frac{1}{\mu^2}.
 \end{align}
\end{thm}

\begin{proof}
 Using Eqs. \eqref{levi} and \eqref{k2}, we have
\begin{align}
\nabla_{\upsilon'}\upsilon'=\left(\phi'+\frac{\sin\phi}{\mu^2}\right)\left(-\frac{\sin\phi}{\mu}\partial_1+\frac{\cos\phi}{\mu}\partial_2+2(\cos\phi)\partial_3\right).
\end{align}
Employing Frenet formula, the above equation yields $\kappa=\phi'+\frac{\sin\phi}{\mu^2}$. If $\kappa>0$ everywhere on $I$ then differentiating normal vector field $N$ along $\upsilon$, we get
\begin{align}\label{quasin}
\nabla_{\upsilon'}N&=-\left(\phi'+\frac{\sin\phi}{\mu^2}\right)\left(\frac{\cos\phi}{\mu}\partial_1+\frac{\sin\phi}{\mu}\partial_2\right)+\left(\frac{\cos2\phi}{\mu^2}-2\phi'\sin\phi\right)\partial_3\nonumber\\
&=-\kappa\upsilon'+\frac{1}{\mu^2}\partial_3.
\end{align}
By the virtue of \eqref{4.1}, we have from Eq. \eqref{quasin} that $\varepsilon\tau B=\frac{1}{\mu^2}\partial_3$, which provides $\tau$. This completes the proof. 
\end{proof}

\section{Null Legendre curves in $\QQ^3$}
If we consider the Reeb vector field to be timelike, i.e. $\varepsilon=-1$ in $\QQ^3$ then the signature of the associated metric \eqref{2.3} is $(2,1)$. In this situation we have lightlike (null) vector fields in $\QQ^3$. Here we give characterization of null Legendre curves in $\QQ^3$. 

\noindent An arbitrary curve $\upsilon:I\to \QQ^3$  is said to be a null curve if
\begin{align}\label{nc}
g(\upsilon'(s),\upsilon'(s))=0\ {\rm and}\ \upsilon'(s)\ne 0\ \forall s\in I.
\end{align}
The general Frenet frame $\{T:=\upsilon', U, V\}$ along null curve $\upsilon$ is determined by
\begin{align}\label{nf}
\begin{cases}
g(U,U)=1,\ g(T,T)=g(V,V)=0,\\ g(T,V)=-1,\ g(T,U)=g(U,V)=0,
\end{cases}
\end{align}
where $T\times U=-T,\, U\times V=-V,\,V\times T=U$. This frame is positively oriented if $det(T,U,V)=1$. The frame  $\{T:=\upsilon', U, V\}$ satisfying the following general Frenet equations: 
\begin{align}\label{nf1}
  \nabla_{T}T =\tilde{h}T+\tilde{ \kappa}_1 U, \,\,\, \nabla_{T}U = -\tilde{\tau}_1  T+\tilde{ \kappa}_1 V \,\,\, {\rm{and}} \,\,\, \nabla_{T}V = -\tilde{h}V-\tilde{\tau}_1 U,
\end{align}
where $\tilde{h},\tilde{ \kappa}_1$ and $\tilde{ \tau}_1$ are smooth functions on $I$. Furthermore, $\tilde{ \kappa}_1$ and $\tilde{ \tau}_1$ are known as curvature and torsion of null curve $\upsilon$, respectively. Since frame $\{T:=\upsilon', U, V\}$ and Frenet formulas \eqref{nf1} depend on the choice of screen distribution and parametrization of $\upsilon$. Therefore, they are not unique in general (for more details see \cite{KLD}).

\begin{prop}
Let $\upsilon:I\to \QQ^3$ be a null Legendre curve in $\QQ^3$. Then $\upsilon$ is a geodesic.
\end{prop}
\begin{proof}
From Eq. \eqref{nf}, we can express 
\begin{align}\label{express}
\nabla_{\upsilon'}\upsilon'= a_3\upsilon'+b_3 U,
\end{align} 
where $a_3$ and $b_3$ are certain smooth functions on $I$. Using Eqs. \eqref{3.1} and \eqref{express} in the differentiation of $g(\xi, \upsilon')=0$ along $\upsilon$, we find that $b_{3}g(\xi, U)=0$. For $b_{3}\neq 0$ we obtain $g(\xi, U)=0$. Then Eq. \eqref{nf} clearly indicates that $\xi=a_{4}\upsilon'+b_{4}V$ along $\upsilon$ for some smooth functions $a_4$ and $b_4$ on $I$. Now using above expression of $\xi$ and Eq. \eqref{nf}, we get $g(\xi, \xi)=-2 a_{4} b_{4}=-1$ which is not possible due to the fact that  $g(\xi, \upsilon')=-b_{4}=0$. Hence, we must have $b_{3}=0$ and therefore \eqref{express} gives $\nabla_{\upsilon'}\upsilon'= a_3\upsilon'$, this implies that $\upsilon$ is a geodesic after a reparametrization. This completes the proof.
\end{proof}

\noindent Below we give the spherical characterization of Legendre curve in a quasi-Sasakian pseudo-metric $3$-manifold $\QQ^3$.
\section{Spherical Legendre curves in $\QQ^3$}
Let $\QQ^3$ be a quasi-Sasakian pseudo-metric $3$-manifold for which $\alpha\ne 0$. Then the sphere $\mathbb{S}_{\pm}^2$ in $\QQ^3$ is defined by $\mathbb{S}_{\pm}^2:=\mathbb{S}_{\circ}^2\cup\mathbb{H}_{\circ}^2$, where $$\mathbb{S}_{\circ}^2=\{P\in M^3\,\vert\, g(P,P)=r^2\}$$ and $$\mathbb{H}_{\circ}^2=\{P\in M^3\,\vert \,g(P,P)=-r^2\}.$$ 
\begin{prop}\label{prop1}
Let $\upsilon:I\to \QQ^3$ be a non-geodesic Legendre curve in $\QQ^3$ such that $\alpha\ne 0$. Then the center of osculating sphere of $\upsilon$ at the point $\upsilon(s)$ is given by 
\begin{align}\label{cent}
c(s)=\upsilon(s)+\frac{1}{\vartheta}N-\frac{\vartheta'}{\vartheta^2\vert\alpha\vert}B.
\end{align}
\end{prop}
\begin{proof}
Let us consider the sphere $\mathbb{S}_{\pm}^2$ passing through four neighbouring points of the Legendre curve $\upsilon$. Let $c(s)$ be the center of sphere and $r$ be its radius. Then the function $\gamma:I\to\RR$ defined by $$\gamma(s)=g(c(s)-\upsilon(s),c(s)-\upsilon(s))\mp r^2$$
satisfies the following conditions:
\begin{align*}
\gamma(s)=\gamma'(s)=\gamma''(s)=\gamma'''(s)=0.
\end{align*}
In light of Eqs. \eqref{4.1} and \eqref{quasikt}, the above equation leads to the following relations:
\begin{align}\label{os}
\begin{cases}
g(c(s)-\upsilon(s),c(s)-\upsilon(s))=\pm r^2,\\
g(c(s)-\upsilon(s),T)=0,\\
g(c(s)-\upsilon(s),N)=\frac{1}{\vartheta},\\
g(c(s)-\upsilon(s),B)=-\frac{\varepsilon\vartheta'}{\vartheta^2\vert\alpha\vert}.
\end{cases}
\end{align}
Since $c(s)-\upsilon(s)\in {\rm Span}\{T,N,B\}$. Therefore, in view of Eq. \eqref{os}, we obtain \eqref{cent}.
\end{proof}
Let the sphere $\mathbb{S}_{\pm}^2$ in $\QQ^3$ be centered at origin and Legendre curve $\upsilon$ lie on it. Then in view of above proposition,  osculating sphere of $\upsilon$ is $\mathbb{S}_{\pm}^2$, $\forall s\in I$.
\begin{prop}\label{prop2}
Let $\upsilon:I\to \QQ^3$ be a non-geodesic Legendre curve in $\QQ^3$ such that $\alpha\ne 0$. Let $\vartheta(s)$ be non-constant if Reeb vector field $\xi$ is spacelike and $\vartheta(s)\ne \vartheta(s_0)\exp\left(\pm\int_{s_0}^{s}\vert\alpha(t)\vert dt\right)$ if $\xi$ is timelike. Then the radius of $\mathbb{S}_{\pm}^2$ is constant, $\forall s\in I$ if and only if the centers of osculating spheres are the same constants. 
\end{prop}
\begin{proof}
Let us assume that the radius of osculating sphere is constant. Then from Eq. \eqref{cent}, we find
\begin{align}\label{os1}
\left(\frac{1}{\vartheta}\right)^2+\varepsilon\left(\frac{\vartheta'}{\vartheta^2\vert\alpha\vert}\right)^2=\pm r^2.
\end{align}
 Differentiating above expression along $\upsilon$, we have
\begin{align*}
\frac{\varepsilon \vartheta'}{\vartheta^{2}\vert\alpha\vert}\left\{\left(\frac{\vartheta'}{\vartheta^{2}\vert\alpha\vert}\right)'-\frac{\varepsilon\vert\alpha\vert}{\vartheta}\right\}=0.
\end{align*}
Consequently, if $\vartheta'\neq 0$ then we get
\begin{align}\label{os2}
\left(\frac{\vartheta'}{\vartheta^{2}\vert\alpha\vert}\right)'-\frac{\varepsilon\vert\alpha\vert}{\vartheta}=0.
\end{align}
Further, differentiation of \eqref{cent} along $\upsilon$ gives
\begin{align}\label{os3}
\nabla_{\upsilon'}c(s)=\left\{\left(\frac{\vartheta'}{\vartheta^{2}\vert\alpha\vert}\right)'-\frac{\varepsilon\vert\alpha\vert}{\vartheta}\right\}B.
\end{align}
In view of Eqs.\eqref{os2} and \eqref{os3}, we find that $c(s)=$ constant for all $s \in I$.\\
Conversely, let $\nabla_{\upsilon'}c(s)=0\  \forall s\in I$ then using \eqref{os3} we obtain
\begin{align}\label{os4}
\left(\frac{\vartheta'}{\vartheta^{2}\vert\alpha\vert}\right)'-\frac{\varepsilon\vert\alpha\vert}{\vartheta}=0.
\end{align}
Now differentiating \eqref{os1} along $\upsilon$, we have
\begin{align*}
\frac{\varepsilon \vartheta'}{\vartheta^{2}\vert\alpha\vert}\left\{\left(\frac{\vartheta'}{\vartheta^{2}\vert\alpha\vert}\right)'-\frac{\varepsilon\vert\alpha\vert}{\vartheta}\right\}=\pm r(s)\left(\nabla_{\upsilon'}r(s)\right).
\end{align*}
Employing \eqref{os4} in the above equation, we get
\begin{align*}
 r(s)\left(\nabla_{\upsilon'}r(s)\right)=0
\end{align*}
which implies that $r(s)=$ constant for all $s\in I$. This completes the proof.
\end{proof}
\noindent By the consequences of Propositions \eqref{prop1} and \eqref{prop2}, we have the following result:
\begin{thm}
Let $\upsilon:I\to \QQ^3$ be a non-geodesic Legendre curve in $\QQ^3$ such that $\alpha\ne 0$. Suppose $\vartheta(s)$ is non-constant if $\xi$ is spacelike and $\vartheta(s)\ne \vartheta(s_0)\exp\left(\pm\int_{s_0}^{s}\vert\alpha(t)\vert dt\right)$ if $\xi$ is timelike. Then the following statements are mutually equivalent: 
\begin{itemize}
\item[(a)] $\upsilon\subset \QQ^3$ is a spherical curve.
\item[(b)] the centers of osculating spheres $\forall s\in I$ at the point $\upsilon(s)$ are the same constants.
\item[(c)] $\left(\frac{\vartheta'}{\vartheta^{2}\vert\alpha\vert}\right)'-\frac{\varepsilon\vert\alpha\vert}{\vartheta}=0$.
\end{itemize}
\end{thm}
\subsection{Integration of Eq. \eqref{os2}}
\noindent Eq. \eqref{os2} can be rewritten as
\begin{align}\label{os5}
\left(\left(\frac{1}{\vartheta}\right)'\frac{1}{\vert\alpha\vert}\right)'+\varepsilon\left(\frac{1}{\vartheta}\right)\vert\alpha\vert=0.
\end{align}
If we consider $\frac{1}{\vert\alpha\vert}=x(s)$ and $\frac{1}{\vartheta}=y(s)$, then above equation takes the form
\begin{align}\label{os6}
\left(x(s)y'(s)\right)'+\varepsilon\left(\frac{y(s)}{x(s)}\right)=0.
\end{align}
Changing the variable in above expression by substituting $\int_{s_0}^{s}\frac{1}{x(t)}dt=z(s)$, then we have
 \begin{align}\label{os7}
\frac{1}{x(s)}\frac{d^2y}{dz^2}=\left(x(s)y'(s)\right)'.
\end{align}
From Eqs. \eqref{os6} and \eqref{os7}, we find
 \begin{align}\label{os8}
\frac{d^2y}{dz^2}+\varepsilon y=0.
\end{align}
If $\xi$ is spacelike then general solution of differential equation \eqref{os8} is given by
\begin{align}\label{os9}
y(s)=A_1\cos\left(\int_{s_0}^{s}\vert\alpha(t)\vert dt\right)+A_2\sin\left(\int_{s_0}^{s}\vert\alpha(t)\vert dt\right).
\end{align}
If $\xi$ is timelike then general solution of \eqref{os8} is given by
\begin{align}\label{os10}
y(s)=B_1\cosh\left(\int_{s_0}^{s}\vert\alpha(t)\vert dt\right)+B_2\sinh\left(\int_{s_0}^{s}\vert\alpha(t)\vert dt\right).
\end{align}
Since $\vartheta(s)\ne \vartheta(s_0)\exp\left(\pm\int_{s_0}^{s}\vert\alpha(t)\vert dt\right)$, therefore we find that $B_1\ne\pm B_2$. 
\begin{rem}
In light of Eqs. \eqref{os9} and \eqref{os10}, it can be seen that Legendre curve given by \eqref{k2} in $\mathcal{Q}^3_{\alpha}$ is not spherical.
\end{rem}



\end{document}